\documentclass[12pt]{amsart}
\usepackage{amssymb,amsmath}
\usepackage{graphicx}
\usepackage{colordvi}

\numberwithin{equation}{section}

\newtheorem{theorem}[equation]{Theorem}

\newtheorem{proposition}[equation]{Proposition}
\newtheorem{cor}[equation]{Corollary}

\theoremstyle{definition}

\newtheorem{example}[equation]{Example}

\copyrightinfo{}{}
\newcounter{FNC}[page]
\def\fauxfootnote#1{{\addtocounter{FNC}{2}$^\fnsymbol{FNC}$%
     \let\thefootnote\relax\footnotetext{$^\fnsymbol{FNC}$#1}}}
\newcommand{\C}{{\mathbb C}}
\renewcommand{\P}{{\mathbb P}}
\newcommand{\PP}{{\mathbb P}}

\newcommand{\Yr}{Y_r}
\newcommand{\CYr}{C\Yr}
\newcommand{\productspace}{\P(S^2V^*)\times_{\C} \P(S^2V)}
\newcommand{\PQU}{\P(S^2Q^*)\times_{G} \P(S^2U)}

\begin{document}

\title{A general formula for the algebraic degree in semidefinite programming}

\author[v. Bothmer]{Hans-Christian Graf von Bothmer}
\address{Institiut f\"ur algebraische Geometrie, Leibnitz Universit\"at Hannover, 
Welfengarten 
1, D-30167 Hannover }
\email{bothmer@math.uni-hannover.de}
\urladdr{http://www.iag.uni-hannover.de/~bothmer/}

\author[K.~Ranestad]{Kristian Ranestad}
\address{Matematisk institutt\\
         Universitetet i Oslo\\
         PO Box 1053, Blindern\\
         NO-0316 Oslo\\
         Norway}
\email{ranestad@math.uio.no}
\urladdr{http://www.math.uio.no/\~{}ranestad}

\subjclass[2000]{14}
\begin{abstract}
 In this note, we use a natural desingularization of the conormal variety of the variety of $(n\times n)$-symmetric matrices of rank at most $r$ 
 to find a general formula for the algebraic degree in semidefinite programming.

\end{abstract}
\maketitle
\section{The algebraic degree in semidefinite programming}

Let $\P^m$ be a general projective space of symmetric $(n\times n)-$matrices  up to scalar multiples, 
and let $Y_{r}\subset \P^m$ be the subvariety of matrices of rank at most $r$.
In this note we find a formula for the degree $\delta(m,n,r)$ of the dual variety $Y_{r}^*$ whenever this is a hypersurface, i.e. 
$\delta(m,n,r)$ is the number of hyperplanes in $\P^m$ in a general pencil that are tangent to $Y_{r}$ at some smooth point.

In \cite {NRS} the algebraic degree of semi-definite  
programming is introduced and is shown to coincide with a degree of a dual variety as above \cite[Theorem 5.2]{NRS}. 
 In particular, the number $\delta(m,n,r)$ is a measure of the algebraic complexity of the rank $r$ solution for a general objective function 
 optimized over a generic $m$-dimensional affine space of symmetric $(n \times n)-$matrices.  

Our general formula for $\delta(m,n,r)$
extends the results of \cite {NRS}.
It is expressed in terms of a function on subsequences of  $\{1,...,n\}$.
 Let 
 $$\psi_{i}=2^{i-1},\quad\psi_{i,j}=\sum_{k=i}^{j-1}\binom{i+j-2}k\quad {\rm when}\quad i<j,$$
 and 
$$\psi_{i_{1},...,i_{r}}={\rm Pf}(\psi_{i_{k},i_{l}})_{1\leq k<l\leq r}\quad \hbox{\rm if $r$ is even}$$
$$\psi_{i_{1},...,i_{r}}={\rm Pf}(\psi_{i_{k},i_{l}})_{0\leq k<l\leq r}\quad \hbox{\rm if $r$ is odd }$$
where $\psi_{i_{0},i_{k}}=\psi_{i_{k}}$ and
Pf denotes the Pfaffian.  
 \begin{theorem}\label{degree}
 The algebraic degree 
$$\delta(m,n,r)=\sum_{I} \psi_{I}\psi_{I^c},$$
where the sum runs over all strictly increasing subsequences $I=\{i_{1},...,i_{n-r}\}$ of $\{1,...,n\}$ of length $n-r$ and sum $i_{1}+...+i_{n-r}=m$, and $I^c$ is the complement 
$\{1,...,n\}\setminus I$.  
\end{theorem}

Notice that the formula makes sense only when $m\geq \binom{n-r+1}2$ and $\binom{r+1}2\leq \binom{n+1}2-m$, 
which is precisely the allowable range for the algebraic degree (cf. \cite[Proposition 3.1]{NRS}).

\section{The bidegree of a conormal variety}

Let $V$ and $V^*$ be an $n$-dimensional vector space and its dual with fixed dual bases.  
Then we interpret a symmetric $(n\times n)-$matrix as a bilinear map $\phi_S:V\to V^*$ (or $\phi_T:V^*\to V$) such that 
$\phi(v_{1})(v_{2})=\phi(v_{2})(v_{1})$ for any pair of vectors $v_{1},v_{2}$. 
 The symmetric maps $V\to V^*$ form $S^2V^*$, while the symmetric maps $V^*\to V$ form $S^2V$.
In this notation we let $\Yr \subset \PP(S^2 V)$ be the space of symmetric $(n \times n)-$matrices of rank at most $r$ (up to scalars) 
and $\CYr \subset \productspace$ the variety of pairs of symmetric 
matrices $(S,T)$ of rank at most $n-r$ and $r$ whose matrix product $ST = 0$. 
The variety $\CYr$ is the conormal variety of $\Yr$, i.e. the closure of the incidence set of pairs $(S,T)$  
such that $S\in \P(S^2V^*)$ corresponds to a hyperplane in $\P(S^2V)$ that is tangent to $Y_{r}$ at the smooth point $T$.

The bidegree of $\CYr\subset \productspace$ is the generating function for the algebraic degree.
More precisely, denote by $H_{1}$ be the class of a hyperplane on $\P(S^2V^*)$ 
and $H_{2}$ be the class of a hyperplane on $\P(S^2V)$, then the class of $\CYr$ is a homogeneous polynomial in $H_{1}$ and $H_{2}$.
\begin{proposition}\cite[Theorem 4.2]{NRS}
  $$[\CYr]=\sum_{m}\delta(m,n,r)H_{1}^{\binom{n+1}2-m}H_{2}^{m}$$ 
\end{proposition}

We may therefore compute the algebraic degree on $\CYr$:
\begin{cor}\label{cor}
\[
\delta(m,n,r) = \int_{\CYr}H_{1}^{m-1}H_{2}^{\binom{n+1}2-m-1}\cap \CYr.
\]
\end{cor}

\section{A desingularization of the conormal variety }
Following the argument of \cite[Lemma 3.2]{NRS} we define a natural desingularization of $\CYr$:
Fix a rank $r$ subspace $K$ of $V$.  The set of pairs $(S,T)$  of symmetric bilinear maps $S:V\to V^*$ and $T:V^*\to V$ 
such that $K\subset \ker S$ and $K^{\bot}\subset \ker T\subset V^*$
form a product of projective spaces:  $S$ determines an element in $S^2(V/K)^*\subset S^2V^*$, 
while $T$ determines an element in $S^2K\subset S^2V$, since $(V^*/K^{\bot})$ is naturally isomorphic to $K^*$.  

Varying $K$ we find a desingularization of $\CYr$ as a fiber bundle over the Grassmannian of rank $r$ subspaces of $V$. 
More precisely, consider the Grassmannian $G=G (r,V)$ of rank $r$ subspaces,
and the universal sequence 
$$0\to U\to V_G\to Q \to 0$$
on $G$.  Consider the symmetric squares $S^2U$ and $S^2Q^*$ and their projectivisations $\P(S^2U)$ and $\P(S^2Q^*)$.  
The inclusions $S^2U\subset S^2V$ and $S^2Q^*\subset S^2V^*$ define a natural map 
$$X:=\PQU\to 
\productspace.$$ 
\begin{proposition} Let $G=G (r,V)$ and $X=\PQU$, where $Q$ and $U$ are the universal quotient and subbundle on $G$.
  Then the image of the natural projection  $X\to 
\productspace$ is the conormal variety $\CYr$ of the variety $\Yr$ of symmetric bilinear maps $T:V^*\to V$ of rank at most $r$, and 
$X\to \CYr$ is a desingularization of $\CYr$.
\end{proposition}
\begin{proof} On the one hand, the image of the projection consists of all pairs $(S,T)$, such that the composition $ST=0$ and the ranks of $S$ and $T$
 are at most $n-r$ and $r$, so the image is precisely the conormal variety $\CYr$ of $\Yr$.  On the other hand  the variety $X$ is clearly smooth, 
and the map $X\to \CYr$ is birational, so it is a desingularization of $\CYr$.
\end{proof}

\section{The proof of the formula}

We denote by $H_{1}$ and $H_{2}$ also  their pullback to $X$. 
It remains by Corollary \ref{cor} to compute $H_{1}^{a}H_{2}^b$ on $X$, 
when $a+b=\binom {n+1}2-2=r(n-r)+\binom {r+1}2-1+\binom {n-r+1}2-1=\dim X$.
 Let $p_{1}: \P(S^2Q^*)\to G$ and $p_{2}: \P(S^2U)\to G$ be the natural projections 
and $p:X\to G$ be the induced projection from the fiber product. 
We use Segre classes to move the computation from $X$ to $G$.   Our projective bundles are the varieties of rank one subbundles and not 
rank one quotient bundles as in \cite{LLT}, so
 the Segre classes in \cite[2.4]{LLT}  are the Segre classes of the dual bundles in our notation.
Since $S^2Q^*$ has rank $\binom {n-r+1}2$, we therefore get
$$ {(p_{1})}_{*} \bigl(H_{1}^{\binom {n-r+1}2-1+i}\cap X \bigr)=(-1)^{i} s_{i}S^2Q^*\cap G= s_{i}S^2Q\cap G.$$

Thus 
$${(p_{1})}_{*}(H_{1}^{a}\cap X) = s_{a-\binom {n-r+1}2+1}S^2Q\cap G.$$
Similarly we compute $H_{2}^b$:
 $${(p_{2})}_{*}(H_{2}^{b}\cap X)=s_{b-\binom {r+1}2+1}S^2U^*\cap G.$$
On the fiber product we get
$$p_{*}(H_{1}^{a}H_{2}^b\cap X)=p_{*}(H_{1}^{a}H_{2}^b\cap p^*G)=s_{a-\binom {n-r+1}2+1}(S^2Q)s_{b-\binom {r+1}2+1}(S^2U^*)\cap G.$$
Combined with Corollary \ref{cor} we have therefore reduced the computation of $\delta(m,n,r)$ to a calculation on 
the Grassmannian variety $G$.

\begin{proposition}The algebraic degree 
$$\delta(m,n,r)=\int_{G}s_{m-\binom {n-r+1}2}(S^2Q)s_{\binom{n+1}2-m-\binom {r+1}2}(S^2U^*)\cap G.$$
\end{proposition}
\begin{proof}   Since $X$ is birational to $\CYr$ and $H_{1}^{m-1}H_{2}^{\binom{n+1}2-m-1}\cap \CYr$ is finite,  
   $$\delta(m,n,r)=\int_{\CYr}H_{1}^{m-1}H_{2}^{\binom{n+1}2-m-1}\cap \CYr=\int_{X} H_{1}^{m-1}H_{2}^{\binom{n+1}2-m-1}\cap X.$$
\end{proof}
In this formulation, computing the degree involves Schubert calculus on $G$.  
In two steps we will circumvent this, following some basic results on Schur functions on Grassmannians (cf. \cite{L} and \cite{LLT}) :
First, we express the total Segre classes of the symmetric product with the Segre classes $s_{i}=s_{i}(E)$ of the original bundle:
$$s(S^2E)=\sum_{I}\psi_{I}s_{I}(E)$$
where $I$ runs over all strictly increasing sequences $\{i_{1},...,i_{e}\}$ of positive integers of length equal to the rank $e$ of $E$,
 and $s_{I}(E)$ is the Schur function
$$s_{I}(E)=\det \begin{pmatrix}
s_{i_{1}-1}&s_{i_{2}-1}&...&s_{i_{e}-1}\\
s_{i_{1}-2}&s_{i_{2}-2}&...&s_{i_{e}-2}\\
...\\
s_{i_{1}-e}&s_{i_{2}-e}&...&s_{i_{e}-e}\\
\end{pmatrix}.$$ 
In particular $s_{I}(E)$ has degree $i_{1}+...+i_{e}-\binom{e+1}2$.
The coefficient $\psi_{I}$ is the integral valued function 
$$\psi_{I}=\sum_{J}A_{I,J},$$  
where $J$ runs over all strictly increasing sequences of length equal to the length of $I$, and
$A_{I,J}$ is the minor of following matrix with rows indexed by $I$ and columns indexed by $J$ 
$$
A=\begin{pmatrix}
1&0&0&0&0&0&...\\
1&1&0&0&0&0&...\\
1&2&1&0&0&0&...\\
1&3&3&1&0&0&...\\
1&4&6&4&1&0&...\\
....&&&&&&...\\
\end{pmatrix}
$$
We refer to \cite[Proposition A.15]{LLT}, for the following closed formulas for the function $\psi_{I}$: 
$$\psi_{i}=2^{i-1},\quad\psi_{i,j}=\sum_{k=i}^{j-1}\binom{i+j-2}k\quad {\rm when}\quad i<j,$$
 and 
$$\psi_{i_{1},...,i_{r}}={\rm Pf}(\psi_{i_{k},i_{l}})_{1\leq k<l\leq r}\quad \hbox{\rm if $r$ is even}$$
$$\psi_{i_{1},...,i_{r}}={\rm Pf}(\psi_{i_{k},i_{l}})_{0\leq k<l\leq r}\quad \hbox{\rm if $r$ is odd }$$
where $\psi_{i_{0},i_{k}}=\psi_{i_{k}}$ and Pf denotes the Pfaffian.

Thus, 
$$s_{m-\binom {n-r+1}2}(S^2Q)=\sum_{I}\psi_{I}s_{I}(Q)$$
where $I$ runs over all strictly increasing sequences of positive integers $I=\{i_{1},i_{2},....,i_{n-r}\}$ with sum $i_{1}+...+i_{n-r}=m$.  
(Notice that our sequences are shifted once to the right compared to those of \cite{LLT}, in particular they do not contain $0$.)

Secondly, if $I$ has length $n-r$ and $J$ has length $r$, then 
$$s_{J}(U)s_{I}(Q)=s_{JI}V_{G}$$
where $JI$ is the concatenated sequence (cf. \cite[3.1]{LLT}).  Now, $s_{JI}=0$ whenever the sequence $JI$ has repeated entries. 
Furthermore, since $V_{G}$ is a trivial bundle,  $s_{JI}V_{G}\cap G=0$ whenever the degree of $s_{JI}V_{G}$ is positive, i.e. the sum of the entries in $JI$ exceeds $\binom{n+1}2$.  
We conclude that  $s_{JI}V_{G}\cap G$ is nonzero only if $JI$ is a permutation of $\{1,...,n\}$ .
In particular, $i_{n-r}, j_{r}\leq n$ and, by (cf. \cite[3.2.1]{LLT})
$$s_{I}(Q)s_{J}(U^*)\cap G=\delta_{I{J^c}}$$
where $I$ and $J$ are subsequences of $\{1,...,n\}$ of length $n-r$ and $r$ respectively, and the sequence $J^c$ is the complement 
$\{1,...,n\}\setminus J$ and, finally, $\delta_{AB}=1$ if $A^c=B$ and $\delta_{AB}=0$ otherwise.

We put the two steps together and the Theorem \ref{degree} immediately follows.

\begin{example}
We compute $\delta(m,n,r)$, when $n=5$.  The relevant values of $\psi$ are:

\[
\begin{array}{cccccc}
I&1&2&3&4&5\\
\psi_{I}&1&2&4&8&16\\
\end{array}
\]

\[
\begin{array}{ccccccccccc}
I&1,2&1,3&1,4&1,5&2,3&2,4&2,5&3,4&3,5&4,5 \\
\psi_{I}&1&3&7&15&3&10&25&10&35&35 \\
\end{array}
\]

\[
\begin{array}{ccccccccccc}
I&1,2,3&1,2,4&1,2,5&1,3,4&1,3,5&1,4,5&2,3,4&2,3,5&2,4,5&3,4,5 \\
\psi_{I}&1&4&11&6&23&27&4&18&30&20 \\
\end{array}
\]

\[
\begin{array}{cccccc}
I&1,2,3,4&1,2,3,5&1,2,4,5&1,3,4,5&2,3,4,5\\
\psi_{I}&1&5&10&10&5\\
\end{array}
\]

So
$$\delta(1,5,4)=\psi_{1}\psi_{2,3,4,5}=5\quad \delta(2,5,4)=\psi_{2}\psi_{1,3,4,5}=20$$
$$\delta(3,5,4)=\psi_{3}\psi_{1,2,4,5}=40\quad \delta(4,5,4)=\psi_{4}\psi_{1,3,4,5}=40$$
$$\delta(5,5,4)=\psi_{5}\psi_{1,2,3,4}=16\quad \delta(3,5,3)=\psi_{1,2}\psi_{3,4,5}=20$$
$$\delta(4,5,3)=\psi_{1,3}\psi_{2,4,5}=90$$
$$\delta(5,5,3)=\psi_{1,4}\psi_{2,3,5}+\psi_{2,3}\psi_{1,4,5}=7\cdot18+3\cdot27=207$$
$$\delta(6,5,3)=\psi_{1,5}\psi_{2,3,4}+\psi_{2,4}\psi_{1,3,5}=15\cdot4+10\cdot23=290$$
$$\delta(7,5,3)=\psi_{2,5}\psi_{1,3,4}+\psi_{3,4}\psi_{1,2,5}=25\cdot6+10\cdot11=260$$
$$\delta(8,5,3)=\psi_{3,5}\psi_{1,2,4}=140\quad\delta(9,5,3)=\psi_{4,5}\psi_{1,2,3}=35$$

\end{example}

\begin{example} With a computer it takes less then one minute to compute
$$\delta(105,20,10) = 167223927145503062075691969268936976274880$$
\end{example}

\bibliography{bibl}

\end{document}